\documentclass[a4paper,12pt]{article}
\usepackage{PhML_2012}
\usepackage{times}
\usepackage{latexsym}
\usepackage{apacite}
\usepackage{graphicx}
\usepackage{amssymb}
\usepackage{epstopdf}

\usepackage[all, 2cell]{xy}
\title{Univalence and Constructive Identity}
\author{Andrei Rodin}
\address{Institute of Philosophy RAS, Russia - University Paris-Diderot, France\\
Volkhonka str. 14, Moscow, 119991, Russia\\
\texttt{rodin@ens.fr}}
\begin{document}

\maketitle
\thispagestyle{empty}

\begin{abstract}
The non-standard identity concept developed in the Homotopy Type theory allows for an alternative analysis of Frege's famous \emph{Venus} example, which explains how empirical evidences justify judgements about identities and accounts for the constructive aspect of such judgements.
\end{abstract}


\section {Introduction}
According to a popular view expressed, in particular, by Frege \cite{Frege:1962}

\begin{quote}
Identity is a relation given to us in such a specific form that it is inconceivable that various kinds of it should occur (p.254)
\footnote{``Die Identitaet ist eine so bestimmt gegebene Beziehung, dass nicht abzusehen ist, wie bei ihr
verschiedene Arten vorkommen k\"onnen.''}
\end{quote}

An alternative view on identity is suggested by the intensional version of Martin-L\"ofÕs Intuitionistic Type theory \cite{Martin-Lof:1984} and its recent interpretation known as the Homotopy Type theory \cite{Warren:2008}, \cite{Awodey&Warren:2009}, \cite{Awodey:2010}. Unlike earlier attempts to diversify the identity relation (such as Geach's theory of ``Relative Identity'' \cite{Geach:1972}, \cite{Deutsch:2002}) the Homotopy Type theory provides new mathematical means for such a diversification. What remains however problematic is whether or not the Homotopy-Type-theoretic account of identity can be taken seriously from a philosophical viewpoint. Although the Homotopy Type theory apparently suggests us a way to conceive of what Frege claims to be inconceivable, namely, an identity concept allowing for different kinds of identity, we still need to assure that this unusual identity concept indeed deserves the name. Since the notion of identity like other basic logical notions is supposed to be applicable outside the pure mathematics it is essential
to check whether or not the Homotopy-Type-theoretic identity meets this requirement.

The main purpose of this short paper is to show that the Homotopy-Type-theoretic identity indeed makes sense in empirical contexts and thus can be possibly used in empirical sciences. The plan of the paper is the following. After introducing the Homotopy-Type-theoretic notion of identity and clarifying the role of Voevodsky's Univalence Axioms in it I turn to FregeÕs classics ``On Sense and Reference'' and use  (informally) the the Homotopy-Type-theoretic identity for  an alternative treatment of Frege's popular \emph{Venus} example. Although this reasoning does not  qualify as a genuine application of Homotopy Type theory in an empirical science it provides a new constructive view on identity and gives an idea of how such an application may look like.

\section{Higher Identity Types and the Univalence Axiom}

Martin-L\"ofÕs  Intuitionistic Type theory with dependent types  \cite{Martin-Lof:1984} involves two different identity relations. First, we have here a notion of \emph{definitional identity} of types (written $A=B$) and of terms belonging to the same given type ($x = y : A$). Corresponding rules assure that the definitional identity is an equivalence relation and that definitionally identical types and terms are mutually interchangeable through substitution in the usual way. Second, we have here a notion of \emph{propositional identity} $Id_{A}(x, y)$ of terms belonging to the same given type $A$.
\footnote{
The original version of this theory involves four different kinds of identity (\cite{Martin-Lof:1984}, page 59). Following Awodey and Warren \cite{Awodey&Warren:2009}, \cite{Awodey:2010} I simplify the original account by deliberately confusing some syntactic and semantical aspects. Then we are left with the following two forms of identity described above in the main text. The version of type theory presented in \cite{Warren:2008} applies the definitional equality also to \emph{contexts}. For simplicity I don't consider the type-theoretic notion of context in this paper.
}
By the \emph{propositions as types} correspondence $Id_{A}(x, y)$ stands both for a proposition saying that terms $x, y$ are identical and for the type of \emph{proofs} (evidences, witnesses) $p_{i}$ proving this proposition. Given two such proofs $p_{i}, p_{j}$  there is a higher identity type $Id_{Id_{A}(x, y)}(p_{i}, p_{j})$ dependent on type $Id_{A}(x, y)$, which is read as a proposition saying that proofs $p_{i}, p_{j}$ are the same. Terms $r_{i}$ of type $Id_{Id_{A}(x, y)}(p_{i}, p_{j})$ are proofs of this latter proposition. The given Type theory allows for higher identity types of all finite orders.

If two terms $x,y$ of given type $A$ are definitionally equal they are interchangeable through substitution and hence also propositionally equal. If the converse is also the case (i.e. if all propositionally equal terms are definitionally equal) the corresponding version of the theory is called \emph{extensional}; otherwise it is called \emph{intensional}. In the extensional theory the difference between the two kinds of identity  is merely formal and arguably redundant but in the intensional theory it is not.

As long as the intensionality property is formulated negatively (as the lack of extensionality) every model of the extensional Type theory also models the corresponding intensional theory. In 1993  Streicher \cite{Streicher:1993} built an intensional model of Martin-L\"ofÕs theory where identity types of the first level (like $Id_{A}(x, y)$) are modeled with abstract groupoids while all higher identity types are reduced in the sense that the propositional identity of terms of higher identity types implies their definitional identity. (Awodey \cite{Awodey:2010} calls this latter property the ``extensionality one dimension up''). More recently Voevodsky observed that the structure of higher identity types can be modeled  by higher \emph{homotopic} groupoids, i.e., groupoids of paths, homotopies between the paths, 2-homotopies between the path homotopies, etc., of some topological space \cite{Voevodsky:2006}. Since a properly defined infinite-dimensional groupoid (aka $\omega$-groupoid) essentially determines a topological space (of which it is the fundamental groupoid) \cite{Grothendieck:1983}, \cite{Voevodsky&Kapranov:1990} one may think of higher groupoids as homotopy groupoids without loss of generality.

Voevodsky's \underline{Univalence Axiom} says that given two types $A, B$ the identity type $Id_{U}(A, B)$ of paths between $A,B$ (seen as terms of type $U$) is equivalent (in $U$) to the type $Eq_{U}(A, B)$ of equivalences between $A$ and $B$, where  $U$ is an  \emph{universe} thought of as the type of \emph{small} types and where the equivalence in $U$ is defined in the usual way as (a term of) the type of inversible functions of the form $A \rightarrow B$. Observe that in order to compare ($n$)-paths with ($n$)-equivalences one needs to rise the dimension and consider $n+1$-maps (cells) between these things. This shows that the Univalence Axiom cannot be satisfied in an $n$-groupoid model for any finite $n$. However since $\omega + 1 = \omega$ in an infinite-dimensional  ($\omega$)-groupoid $\omega$-equivalences between objects compare themselves with $\omega$-paths as the Univalence requires. Thus the Univalence Axiom rules out all groupoid models of  Martin-L\"ofÕs theory except the infinite-dimensional ($\omega $-groupoid) ones, which can be ``intensional all the way up''.

 \section{Frege and Constructive Identity}
Frege \cite{Frege:1892}, \cite{Frege:1952} develops his theory of meaning (sense) and reference as a response to the following problem:

\begin{quote}
$a = a$ and $a = b$ are obviously statements of different cognitive value; $a = a$ holds \emph{a priori} and, according to Kant, is to be labeled analytic, while statements of the form $a = b$ often contain very valuable extensions of our knowledge and cannot always be established \emph{a priori}. The discovery that the rising sun is not new every morning, but always the same, was one of the most fertile astronomical discoveries. Even today the identification of a small planet or a comet is not always a matter of course. Now if we were
to regard identity as a relation between that which the names $a$ and $b$ designate, it would seem
that $a = b$ could not differ from $a = a$ (provided $a = b$ is true). (\cite{Frege:1952}, p. 56, the English translation is slightly modified)
\end{quote}

After presenting his theory Frege suggests this solution of the above problem:

\begin{quote}
When we found $a = a$ and $a = b$ to have different cognitive values, the explanation is that for the
purpose of knowledge, the sense of the sentence, viz., the thought expressed by it, is no less
relevant than its reference, i.e. its truth value. If now $a = b$, then indeed the reference of $b$ is the
same as that of $a$, and hence the truth-value of $a = b$ is the same as that of $a = a$.  In spite of this the sense of $b$ may differ from that of $a$ and thereby the thought expressed in $a = b$ differs from that of $a = a$. In that case the two
sentences do not have the same cognitive value. If we understand by judgment the advance from
the thought to its truth value, as in the above paper, we can also say that the judgments are different. (\cite{Frege:1952}, p. 78)
\end{quote}

Frege's theory indeed explains how judgements $\vdash a = b$ and $\vdash a = a$ can be both true and still have different cognitive (and epistemic) values. However this theory has two essential shortcomings: (i) it does not specially account for the case when in $\vdash a = b$ terms $a$ and $b$ have the same meaning (expressed by different symbols) and (ii) it does not explain how empirical or other evidences may possibly justify judgement $\vdash a = b$ when terms $a, b$ have different meanings. The latter shortcoming has to do with Frege's adherence to the old-fashioned pre-Kantian view according to which basic logical concepts like that of truth and identity must be fixed independently of the practice of empirical and mathematical research, so the logic of judgements of the form $ \vdash a = b$  is not relevant to the question ``How we come to know that $ a = b$?''. In what follows I suggest an alternative solution of Frege's problem that does take this question into consideration. Of many Frege's examples found in his \cite{Frege:1892}, \cite{Frege:1952} I shall use the most popular \emph{Venus} example:  here $a$ stands for \emph{Morning Star},  $b$ stands for \emph{Evening Star} and both expressions refer to the same planet Venus.

In fact Frege's examples may be straightforwardly analyzed in terms of the aforementioned distinction between definitional and propositional identities in Martin-L\"of's Type theory. The identities  \emph{Morning Star = MS} where $MS$ is an abbreviation of \emph{Morning Star}, and \emph{Evening Star = ES} where $ES$ is an abbreviation of \emph{Evening Star}, are clearly \underline{definitional}, while the identity \emph{Morning Star = Evening Star} is propositional and requires a proof (evidence, reason). Using the Type-theoretic formalism one can write this latter proposition as $Id_{A}(MS, ES)$ (where $A$ is an appropriate type of observable celestial bodies) and also interpret this formula through the \emph{propositions as types} correspondence as the type of evidences that the given identity holds. 

We suppose after Frege that identities of $MS$ and $ES$ are fixed beforehand. This means that we are in a position to identify two independent observations of $MS$  as observations of \emph{one and the same} object without referring to any explicit evidence or reason for this identification; similarly for ES. This assumption marks our chosen starting point: a deeper analysis would require making such evidences/reasons explicit. Now in order to establish $Id_{A}(MS, ES)$ we \underline{construct} invertible maps $m_{i}$ from available observations of $MS$ to available observations of $ES$ (by mapping the corresponding frames of reference). These maps (invertible transformations) allow for discerning  certain invariant properties, which we predicate to a new object that we call Venus; we claim that  (non-relational) properties of Venus are the same when we observe it in the Morning and in the Evening.

Thinking about physical objects in terms of invariants of groups of transformations between different representations of these objects is well-established in the 20th century physics: for an easy example think of the shape of an extended event in the flat relativistic spacetime. However the Homotopy Type theory allows for a more advanced way of object-building. Instead of being satisfied with the fact that at least one appropriate map $m$ of type $Id_{A}(MS, ES)$ exists we study (the space of) such maps more attentively. As long as maps $m_{i}$ count as proofs (evidences) of $Id_{A}(MS, ES)$ we wish to classify and evaluate such evidences. In particular, we need some criteria by which maps $m_{i}, m_{j}$ count as essentially the same or count as different. For this purpose we construct appropriate second-order maps $m^{ij}_{k}$ of types  $Id_{Id_{A}(MS, ES)}(m_{i}, m_{j})$ satisfying the usual coherence conditions and consider them as constitutive elements of  object Venus on equal footing with maps $m_{i}$. Thinking about these and further higher-order maps (of higher identity types) as (higher) homotopies allows us to associate with Venus a topological structure, which reflects the way in which this object is empirically construed from available empirical data. Beware that the construction of a higher-dimensional fundamental groupoid that has been just described is all the way \emph{empirical} in the sense that the coherence conditions under the suggested interpretation are checkable against the observational data.

Returning back to Frege's seminal paper \emph{On Sense and Reference} \cite{Frege:1892}, \cite{Frege:1952} I would like to note that although the author opens this paper with a problem relevant to the astronomical research of his time, his suggested solution is motivated mostly by a logical analysis of historical narrative and everyday talk (such examples, which I don't quote here, are abundant in Frege's paper). In my view this is a wrong strategy because the concept of identity that we need for doing science can be simply not found in linguistic examples of this sort. New approaches to identity emerging in the pure mathematics may, on the contrary, turn to be fruitful in empirical sciences - just like new approaches to the notion of space in the 19th century geometry has turned to be highly fruitful in the 20th century physics. 
I see it as a philosopher's task to discern and develop these new approaches beyond their purely technical aspects.  Although the above analysis of Frege's \emph{Venus} example is a dummy it shows how the non-standard Homotopy-Type-theoretic identity can be thought of in empirical terms; I hope that this informal analysis may facilitate some real applications in the future.

\section{To Weaken or To Construct?}
I would like to conclude by comparing the constructive approach to identity outlined above with another non-standard approach to identity, which is  motivated by the philosophy of Structuralism \cite{Keranen:2001}. This latter approach, as exemplified by KrauseÕs theory of quasi-sets \cite{Krause:1992} and its tentative application in particle physics \cite{French&Krause:2006}, seeks to ``weaken'' the standard Fregean identity by replacing it with another equivalence relation. Some people tend to interpret the Homotopy-Type-theoretic account of identity in the same vein pointing to the fact that the homotopy equivalence is a weaker relation than the strict identity. In my view such an interpretation is rather misleading because it does not take into account the constructive aspect of Martin-L\"of's Type theory. The idea of weak identity  suggests that some part of the content of the standard (``strict'') identity concept is taken out while the rest remains untouched. But this is not how the non-standard homotopic identity is really construed! Observe that in  Martin-L\"of's theory the standard identity relation is accounted for in the form of \emph{definitional} identity; the non-standard \emph{propositional} identity is a further construction on the top of the definitional identity that adds a new conceptual content.  Technically speaking, Krause's non-standard identity relation is also construed on the top of the standard identity, so what I want now to stress is not a specific technical feature but rather the underlying philosophy that guides technical developments. While the structuralist idea of weak identity does not explain why in order to get to the weak identity one needs to begin with the strict identity anyway (which makes one to suspect that the ``real'' identity is the latter rather than the former) the constructive approach to identity does this as follows: we need the standard ``strict'' identity for fixing terms of our reasoning, which is a necessary condition for building further non-trivial identities like \emph{Morning Star = Evening Star} or $7 + 5 = 12$ marking a real advance of our knowledge.

\bibliographystyle{apacite}
\bibliography{univalence_extended_styled}
\end{document}